\documentclass[]{article}
\usepackage{graphicx}
\usepackage{amsthm,amsmath,amssymb,amsfonts,latexsym}

\usepackage{natbib}

\newtheorem{theorem}{Theorem}[section]
\newtheorem{lemma}[theorem]{Lemma}
\newtheorem{corollary}[theorem]{Corollary}
\newtheorem{proposition}[theorem]{Proposition}

\newtheorem{definition}[theorem]{Definition}

\newtheorem{eg}[theorem]{Example}

\newcommand{\hf}[1][1]{\frac{#1}{2}}

\newcommand{\di}{\,\mathrm{d}}
\newcommand{\rmd}{\,\mathrm{d}}
\newcommand{\Or}{\mathrm{O}}
\newcommand{\rmi}{\mathrm{i}}

\begin{document}


\setlength{\unitlength}{10mm}
\title{On the stochastic Burgers equation with some applications to turbulence and astrophysics}
\author{A D Neate \and  A Truman}

\maketitle
\begin{abstract}
We summarise a selection of results on the inviscid limit of the stochastic Burgers equation emphasising geometric properties of the caustic, Maxwell set and Hamilton-Jacobi level surfaces and relating these results to a discussion of stochastic turbulence. We show that for small viscosities there exists a vortex filament structure near to the Maxwell set. We  discuss how this vorticity is directly related to the adhesion model for the evolution of the early universe and include new explicit formulas for the distribution of mass within the shock.

\end{abstract}

\section{Introduction}
The Burgers equation was first introduced by J. M. Burgers as a model for pressureless gas dynamics.  It has since provided a tool for studying turbulence in fluids \citep{MR1999766}, for obtaining detailed asymptotics for stochastic Schr\"{o}dinger and heat equations \citep{MR1370174, MR1401222,ETZ} and has played a part in Arnol'd's work on caustics \citep{MR1345386, MR1151185, MR1178935} and Maslov's works in semiclassical quantum mechanics  \citep{MR634377}. It has also been used for studying the formation of the early universe in the Zeldovich approximation and also the adhesion model \citep{MR657014, MR989562}. A detailed explanation of these applications as well as a complete history of the Burgers equation can be found in \citep{Bec}.

In this article we will summarise a selection of results on the inviscid limit of the stochastic Burgers equation and outline some applications of these results to turbulence and the adhesion model. We begin in Section 2 with a summary of results on deterministic Hamilton-Jacobi theory for the heat and Burgers equation.

In Sections 3 to 5 we present some
 geometric and analytic results first developed by Davies, Truman and Zhao \citep{MR1902482, MR2131267} and later extended by Truman and Neate \citep{MR2167533, MR2314842}.  These results relate the geometry of the caustic, Hamilton-Jacobi level surfaces and Maxwell set to that of their algebraic pre-images under the inviscid classical mechanical flow map $\Phi_t$ which will be defined in Section 3.
In two dimensions these results show that a Hamilton-Jacobi level surface, or Maxwell set can only have a cusp  where their pre-images intersect the pre-caustic and so can only have cusps on the caustic. They also allow us to give conditions for the formation of swallowtails on both caustics and level surfaces which in turn have implications for the geometry of the Maxwell set.

We also introduce a reduced (one dimensional) action function which was developed by Reynolds, Truman and Williams \citep{MR2083376} under the assumption that only singularities of $A_k$ type occur \citep{MR1178935}. Using this, we can find explicit equations for the caustic, level surfaces and Maxwell set and their pre-images. In Section 6 we use this to write down an explicit stochastic process whose zeros give `turbulent times' at which cusps on the Hamilton-Jacobi level surfaces appear and disappear infinitely rapidly.

Finally, in Sections 7 and 8, we summarise results   showing that the fluid has  non-zero vorticity in some neighbourhood of the Maxwell set \citep{MR2314842} . We show that this vorticity disappears under the assumptions required for the adhesion model for the evolution of the early universe and outline a new formula for the mass which adheres to the shock (the Maxwell set).

\noindent\emph{Notation:} Throughout this paper
$x,x_0,x_t$ etc will denote vectors (usually in $\mathbb{R}^d$).
Cartesian coordinates of these will be indicated using a sub/superscript where
relevant; thus $x=(x_1,x_2,\ldots,x_d)$,
$x_0=(x_0^1,x_0^2,\ldots,x_0^d)$ etc. The only exception  will be
 in discussions of explicit examples in two and three dimensions when
we will use $(x,y)$ and $(x_0,y_0)$ etc to denote the vectors.

\section{Elements of Hamilton-Jacobi Theory}
We begin by considering a deterministic classical  mechanical system consisting of a unit mass moving under the influence of a conservative force, $-\nabla V$. This system has Hamiltonian,
\[H(q,p) = \hf p^2 +V(q),\]
where $p,q\in\mathbb{R}^d$. 
 Let us assume that the system has a  given initial velocity field $\nabla S_0$ for some function $S_0:\mathbb{R}^d \rightarrow \mathbb{R}$.

The evolution of this system will be given by the classical mechanical flow map, $\Phi_s:\mathbb{R}^d\rightarrow\mathbb{R}^d$ defined by,
\[\frac{\di^2 \Phi_s}{\di s^2} = -\nabla V(\Phi_s),\]
with initial condition,
\[\Phi_0 = I_d, \qquad \dot{\Phi}_0 = \nabla S_0,\]
where $I_d$ denotes the $d$-dimensional identity map. 
Thus, if $X(s)$ is a classical mechanical path  with $X(0) = x_0$, then,
\[X(s) = \Phi_s (x_0), \quad  \dot{X}(0) = \nabla S_0(x_0).\] Usually we also demand that $X(t) = x$ for fixed $x$ and $t$.
If $S_0$ and $V$ are  twice continuously differentiable with bounded second order derivatives, then there exists a caustic time $t_c>0$, such that for all  
$t\in(0,t_c)$ 
the classical mechanical flow map is a diffeomorphism. This is a simple consequence of the global inverse function theorem \citep{MR515141}.  Therefore we can define,
\[x_0(x,t): = \Phi_t^{-1}(x),\]
to be the unique pre-image of the point $x$ reached by the path $X(s)$ at time $t$.
If we now define,
\[\mathcal{S}(x,t) := S_0(x_0(x,t)) + \int_0^t \left(\hf \dot{X}^2(s) - V(X(s))\right)\di s,\]
then it can be easily shown that $\mathcal{S}_t(x) := \mathcal{S}(x,t)$ satisfies the Hamilton-Jacobi equation,
\begin{equation}\label{HJ1}\frac{\partial \mathcal{S}_t}{\partial t} +H(x,\nabla \mathcal{S}_t)
=0,\quad\mathcal{S}_{t=0}(x) = S_0(x).
\end{equation}
We now show how the function $\mathcal{S}_t$ can be used to construct a semi-classical solution to a corresponding heat equation \citep{MR0468817, MR1370174, MR1628310}.

Consider the heat equation for $u^{\mu}(x,t) \in\mathbb{R}$ where $ x\in\mathbb{R}^d$ and $t>0$,
\begin{equation}\label{HE1}\frac{\partial u^{\mu}}{\partial
t}=\hf[\mu^2]\Delta u^{\mu} +\mu^{-2}V(x)u^{\mu},
\end{equation}
with  initial condition,
\begin{equation}\label{HE2}u^{\mu}(x,0) = \exp\left(-\frac{S_0(x)}{\mu^2}\right)T_0(x).
\end{equation}

Let $B_s\in\mathbb{R}^d$ be a $d$-dimensional  Wiener process on the space $(\Omega,\mathcal{F},\mathbb{P})$ with $\mathbb{E}\left\{B(s) B(t) \right\}= \min(s,t)$. Define an Ito diffusion $X_s^{\mu}\in\mathbb{R}^d$ and an Ito process $Y_s^{\mu}\in\mathbb{R}^d$ by,
\begin{alignat}{4}
\di X_s^{\mu}& =  - \nabla \mathcal{S}_{t-s}(X^{\mu}_s)  \di s + \mu \di B_s,\qquad &  X^{\mu}_0 &= x, \label{ItoX}\\
\di Y_s^{\mu}& =  \mu \di B_s, & Y^{\mu}_0 &= x,\label{ItoY}
\end{alignat}
where $0<s\le t<t_c$.  The time reversal in $\mathcal{S}_{t-s}$ allows us to effectively consider a diffusion process which will reach the point $x$ at time $t$. Define $h(s,\omega) := h_0(Y_s^{\mu}(\omega),s)$ where,
\[h_0(Y_s^{\mu},s): = -\mu^{-1}\nabla \mathcal{S}_{t-s}(Y_s^{\mu}).\]
 Since $h$ satisfies the Novikov condition,
\[\mathbb{E}_{\mathbb{P}}\left\{ \exp\left(\hf \int_0^{t_c} h^2(s,\omega)\di s\right)\right\}<\infty,\]
where $\mathbb{E}_{\mathbb{P}}$ denotes expectation with respect to the measure $\mathbb{P}$, it follows that,
\[M_s: = \exp\left(-\int_0^s h(u,\omega) \di B_u - \hf\int_0^sh^2(u,\omega) \di u\right),\]
is a martingale with respect to $\mathcal{F}_s=\sigma(B_s)$ and $\mathbb{P}$.

Using the Girsanov theorem, we can now define a new measure $\tilde{\mathbb{P}}$ on $(\Omega,\mathcal{F})$,
\[\di \tilde{\mathbb{P}} (\omega) = M_{t_c}(\omega) \di \mathbb{P}(\omega),\]
and then,
\[\tilde{B}_s := \int_0^s h(u,\omega) \di u + B_s,\]
is a Brownian motion with respect to $\tilde{\mathbb{P}}$. Therefore, $(Y^{\mu}_s,\tilde{B}_s)$, where $Y_s^{\mu}$ is defined in (\ref{ItoY}),  forms a weak solution to equation (\ref{ItoX}). That is,
\[\di Y_s^{\mu} = - \nabla \mathcal{S}_{t-s}(Y^{\mu}_s) \di s + \mu\di\tilde{B}_s, \]
and conseqeuntly,
\begin{equation}\label{girs}
\mathbb{E}_{\mathbb{P}}\left\{f(X^{\mu}_s)\right\}=\mathbb{E}_{\tilde{\mathbb{P}}}\left\{f(Y^{\mu}_s)\right\}=\mathbb{E}_{\mathbb{P}}\left\{M_s f(B_s)\right\}.
\end{equation}

It follows from the Feynmann-Kac formula that the heat equation (\ref{HE1}) has a solution given by,
\[u^{\mu}(x,t) = \mathbb{E}_{\tilde{\mathbb{P}}}\left\{T_0(Y_t^{\mu}) \exp\left(-\mu^{-2}S_0(Y_t^{\mu}) +\mu^{-2}\int_0^t V(Y_s^{\mu}) \di s\right)\right\},\]
and so by equation (\ref{girs}),
\begin{eqnarray}
u^{\mu}(x,t) &=& \mathbb{E}_{\mathbb{P}}\left\{T_0(X_t^{\mu}) \exp\left(-\mu^{-2}S_0(X_t^{\mu}) +\mu^{-2}\int_0^t V(X_s^{\mu}) \di s\right)\frac{\di \tilde{\mathbb{P}}}{\di\mathbb{P}}\right\}\nonumber\\
& = &  \mathbb{E}_{\mathbb{P}}\left\{T_0(X_t^{\mu}) \exp\left(-\mu^{-2}S_0(X_t^{\mu}) +\mu^{-2}\int_0^t V(X_s^{\mu}) \di s\right.\right.\nonumber\\
& &\quad\left.\left. +\mu^{-1}\int_0^t \nabla \mathcal{S}_{t-s}(X_s^{\mu}) \di B_s - \frac{1}{2\mu^2}\int_0^t|\nabla \mathcal{S}_{t-s}(X_s^{\mu})|^2 \di s\right)\right\}_.\label{FK1}
\end{eqnarray}
Now, using Ito's formula,
\begin{eqnarray*}
 \mathcal{S}(X^{\mu}_t,0) & = & \mathcal{S}(x,t) + \int_0^t \left(\frac{\partial \mathcal{S}_{t-s}}{\partial s}(X^{\mu}_s) -| \nabla \mathcal{S}_{t-s}(X_s^{\mu})|^2 +\hf[\mu^2]\Delta \mathcal{S}_{t-s} (X^{\mu}_s)\right) \di s\\ & &\quad
+\mu \int_0^t \nabla \mathcal{S}_{t-s}(X^{\mu}_s) \di B_s,
\end{eqnarray*}
and so substituting into equation (\ref{FK1})  for $\int_0^t \nabla\mathcal{S}_{t-s}(X^{\mu}_s)\di B_s$ gives,
\begin{eqnarray*}u^{\mu}(x,t) 
& = & \mathrm{e}^{-\frac{\mathcal{S}_t(x)}{\mu^2}} \mathbb{E}_{\mathbb{P}}\left\{T_0(X^{\mu}_t) \exp\left(-\hf\int_0^t \Delta\mathcal{S}_{t-s}(X^{\mu}_s)\di s\right.\right.\\
& &\quad\left.\left. -\mu^{-2}\int_0^t \left(\frac{\partial \mathcal{S}_{t-s}}{\partial s}(X^{\mu}_s) - \hf|\nabla\mathcal{S}_{t-s}(X^{\mu}_s)|^2 - V(X^{\mu}_s)\right) \di s\right)\right\}.
 \end{eqnarray*}
But $\mathcal{S}_t$ satisfies the Hamilton-Jacobi equation (\ref{HJ1}), and so, by reversing time in the diffusion $X^{\mu}$, we have,
\begin{equation}\label{eqn1} u^{\mu}(x,t) = \exp
\left(-\frac{\mathcal{S}_t(x)}{\mu^2}\right) \mathbb{E}_x \left\{T_0(X_0^\mu)
\exp\left(-\hf\int_0^t\Delta \mathcal{S}_{t-s}(X_s^{\mu})\di s\right)\right\}.\end{equation}

Using the logarithmic Hopf-Cole transformation \citep{MR0047234},
\begin{equation}\label{HC1}v^{\mu}(x,t)=-\mu^2\nabla\ln
u^{\mu}(x,t),
\end{equation}
 the heat equation (\ref{HE1}) becomes the Burgers equation for velocity field $v^{\mu}(x,t)\in\mathbb{R}^d$ where $\mu^2$ is now the coefficient of viscosity,
\begin{equation}\label{BE1}\frac{Dv^{\mu}}{D t} = \frac{\partial v^{\mu}}{\partial
t}+(v^{\mu}\cdot \nabla)v^{\mu} =\hf[\mu^2]\Delta v^{\mu} -\nabla V,
\end{equation}
with initial condition,
$$v^{\mu}(x,0) = \nabla S_0(x)+\mathrm{O}(\mu^2).$$

We will be particularly interested in the behaviour of $v^{\mu}$ for small values of $\mu$. In the remainder of this paper we will focus on the discontinuities that develop in $v^{\mu}$ as $\mu\rightarrow 0$.

The convergence factor $T_0$ in the initial condition (\ref{HE2}) is related to the square root of the Burgers fluid mass
density $\rho_t^{\hf}$,
\begin{equation}\label{eqn3}T_0(x_0(x,t))\left|\left(\frac{\partial x_0}{\partial
x}(x,t)\right)\right|^{\hf} = \rho_t^{\hf}(x).\end{equation} 
For $t\in(0,t_c)$ it can be seen that mass is conserved,
$$\mbox{total mass} = \int \rho_t(x)\di x = \int T_0^2(x_0)\di x_0 = \int
\rho_0(x)\di x.$$

The next lemma will be key to our treatment of the solution for the Burgers equation.
\begin{lemma}\label{key1}
Consider the above $C^2$ Hamiltonian dynamical system with Hamiltonian $H(q,p)$ and Hamilton-Jacobi function $\mathcal{S}_t$ satisfying,
\[\frac{\partial \mathcal{S}_t}{\partial t} +H(x,\nabla\mathcal{S}_t) = 0,\qquad \mathcal{S}_{t=0}(x)=S_0(x),\]
so that,
\[\dot{X}(t) = \nabla \mathcal{S}_t(X(t)),\qquad \dot{X}(0) = \nabla S_0(X(0)).\]
Then,
\[\exp\left\{-\hf \int_0^t \Delta \mathcal{S}_s(X(s)) \di s\right\} = \left|\frac{\partial X(0)}{\partial X(t)}\right|^{\hf},\]
where the right hand side is a Jacobian determinant.
\end{lemma}

In particular it follows from Lemma \ref{key1}, for $t\in(0,t_c)$,  that by considering an asymptotic expansion of the diffusion $X_s^{\mu}$ in the solution to the heat equation (\ref{eqn1}),
\begin{equation}\label{eqn2}u^{\mu}(x,t) = \exp
\left(-\frac{\mathcal{S}_t(x)}{\mu^2}\right)T_0(x_0(x,t))\times\left|\left(\frac{\partial x_0}{\partial
x}(x,t)\right)\right|^{\hf}(1+\mathrm{O}(\mu^2)),\end{equation}
where $x_0(x,t)$ is the unique start point of $X^0_s$ with,
$$\dot{X}^0_s=\nabla\mathcal{S}_s(X_s^0),\qquad X_t^0=x.$$
Consequently, the Burgers velocity field is given by,
$$v^{\mu}=v^{\mu}(x,t) \sim \nabla \mathcal{S}_t(x)+\mathrm{O}(\mu^2).$$

\section{The stochastic case}
We now consider the behaviour of a Burgers equation with stochastic forcing. That is for  $v^{\mu}(x,t)\in\mathbb{R}^d$,
\begin{equation}\label{SBE1}\frac{\partial v^{\mu}}{\partial
t}+\left(v^{\mu}\cdot \nabla\right) v^{\mu} = \hf[\mu^2]\Delta
v^{\mu} - \nabla V(x) -\epsilon \nabla k_t(x)\dot{W}_t,
\end{equation} with
initial condition  $v^{\mu}(x,0) = \nabla S_0 (x)+\Or(\mu^2)$, where $\dot{W}_t$ denotes
white noise.

 Using the logarithmic
Hopf-Cole transformation (\ref{HC1}), the Burgers equation (\ref{SBE1}) becomes the Stratonovich heat
equation,
\begin{equation}\label{SHE1}\frac{\partial u^{\mu}}{\partial
t}=\hf[\mu^2]\Delta u^{\mu} +\mu^{-2}V(x)u^{\mu}+\frac{\epsilon}{\mu^{2}}
k_t(x)u^{\mu}\circ\dot{W}_t ,
\end{equation}
with initial condition  $u^{\mu}(x,0) = \exp\left(-\frac{S_0(x)}{\mu^2}\right)T_0(x).$

Now let,
\[
A[X]: = \hf\int_0^t \dot{X}^2(s)\di s -\int_0^t V(X(s)) \di s
-\epsilon \int_0^t k_s(X(s))\di W_s,\] and select a path $X$ with
$X(t)=x$ which minimises $A[X]$. This requires,
\[
\di \dot{X}(s)+\nabla V(X(s)) \di s
+\epsilon \nabla k_s(X(s))\di W_s=0.
\]
We then define the stochastic action,
$A(X(0),x,t):=\inf\limits_X\left\{ A[X]:X(t)=x\right\}.$
Setting,
\[\mathcal{A}(X(0),x,t):=S_0(X(0))+A(X(0),x,t),\]
and then minimising $\mathcal{A}$  over $X(0)$, gives
$\dot{X}(0) = \nabla
S_0(X(0)).$
Moreover, it follows that, \[\mathcal{S}_t(x):=\inf\limits_{X(0)}
\left\{\mathcal{A}(X(0),x,t)\right\},\]
is  the minimal solution of the Hamilton-Jacobi
equation,
\[
\di \mathcal{S}_t +\left(\hf|\nabla \mathcal{S}_t|^2+V(x)\right)\di
t +\epsilon k_t(x) \di W_t=0, \qquad \mathcal{S}_{t=0}(x) =
S_0(x).\] Following the work of Freidlin and Wentzell
\citep{MR1652127}, \[-\mu^2\ln u^{\mu}(x,t) \rightarrow
\mathcal{S}_t(x),\] as $\mu\rightarrow 0$. This gives the inviscid
limit of the minimal entropy solution of the  Burgers equation  as
$v^0(x,t)=\nabla\mathcal{S}_t(x)$ \citep{MR2169977}.

Define the classical flow map
$\Phi_s:\mathbb{R}^d\rightarrow\mathbb{R}^d$ by,
\[\rmd \dot{\Phi}_s +\nabla V(\Phi_s) \rmd s+\epsilon\nabla k_s(\Phi_s)\rmd
W_s =0,\qquad\Phi_0 = \mbox{id},\qquad \dot{\Phi}_0 = \nabla S_0.\]
Since $X(t) =x$ it follows that $X(s) = \Phi_s\left(
\Phi_t^{-1}( x)\right)\!,$ where the pre-image $x_0(x,t) = \Phi_t^{-1} (x)$ is
not necessarily unique.

Given some regularity and boundedness, the global inverse function
theorem gives a random caustic time $t_c(\omega)$ such that for
$0<t<t_c(\omega)$, the pre-image, $x_0(x,t)$, if it exists, is unique and $\Phi_t$ is a random diffeomorphism. Thus, before the
caustic time $v^0(x,t) = \dot{\Phi}_t\left(\Phi_t^{-1}(x)\right)$ is
the inviscid limit of a solution of the Burgers equation
with probability one \citep{MR1370174, MR1628310}.

The method of characteristics suggests that discontinuities in
$v^0(x,t)$ are associated with the non-uniqueness of the real pre-image $x_0(x,t)$.
In the sitution we consider, when this occurs the classical flow map
$\Phi_t$ focusses
an infinitesimal volume of points $\rmd
x_0$ into a zero volume $\rmd X(t)$.
 \begin{definition}\label{i2} The caustic at time $t$ is defined to be the set,
\[ C_t = \left\{ x: \quad\det\left(\frac{\partial X(t)}{\partial x_0}\right) = 0\right\}. \]
\end{definition}

 Assume that after the caustic time $t_c(\omega)>0$, $x$ has $n$ real pre-images,
\[\Phi_t^{-1}\left\{x\right\} =
\left\{x_0(1)(x,t),x_0(2)(x,t),\ldots,x_0(n)(x,t)\right\},\]  where
each $x_0(i)(x,t)\in\mathbb{R}^d$. Then the Feynman-Kac formula and
Laplace's method in infinite dimensions give for a non-degenerate
critical point \citep{MR692299, MR771765},
 \begin{equation}\label{useries}u^{\mu}(x,t)=
\sum\limits_{i=1}^n \theta_i
\exp\left(-\frac{S_0^i(x,t)}{\mu^2}\right),
\end{equation} where
$S_0^i(x,t) :=
S_0\left(x_0(i)(x,t)\right)+A\left(x_0(i)(x,t),x,t\right),$
and $\theta_i$ is an asymptotic series in $\mu^2$. An asymptotic
series in $\mu^2$  can also be found for
$v^{\mu}(x,t)$ \citep{MR1628310}.
Note that $\mathcal{S}_t(x) = \min
\{S_0^i(x,t):i=1,2,\ldots,n\}$.

\begin{definition}\label{i3}
The Hamilton-Jacobi level surface is the set,
\[H_t^c = \left\{ x:\quad S_0^i(x,t) =c \mbox{ for
some }i\right\}.\]
\end{definition}
As $\mu\rightarrow 0$, the dominant term in the expansion (\ref{useries}) comes
from the minimising $x_0(i)(x,t)$ which we denote $\tilde{x}_0(x,t)$.  Assuming $\tilde{x}_0(x,t)$ is unique, we obtain the inviscid limit of
the Burgers  fluid velocity as the minimal entropy
$v^0(x,t) = \dot{\Phi}_t\left(\tilde{x}_0(x,t)\right).$

If the minimising pre-image $\tilde{x}_0(x,t)$ suddenly changes
value between two pre-images $x_0(i)(x,t)$ and $x_0(j)(x,t)$, a jump
discontinuity will occur in $v^0(x,t)$. There are  two distinct ways
in which the minimiser can change; either two pre-images coalesce
and disappear (become complex), or the minimiser switches between
two pre-images at the same action value. The first of these occurs
as $x$ crosses the caustic. When this results in the minimiser
disappearing the caustic is said to be cool. The second occurs as
$x$ crosses the Maxwell set and again, when the minimiser is
involved, the Maxwell set is said to be cool.

 \begin{definition}\label{m1}
The Maxwell set is,
\begin{eqnarray*}
M_t & = & \left\{x:\, \exists \,x_0,\check{x}_0\in\mathbb{R}^d \mbox{ s.t. }\right.\\
&   &\quad\left.
x=\Phi_t(x_0)=\Phi_t(\check{x}_0), \,x_0\neq \check{x}_0 \mbox{ and }
 \mathcal{A}(x_0,x,t)=\mathcal{A}(\check{x}_0,x,t)
\right\}.
\end{eqnarray*}
 \end{definition}

\begin{figure}[h]
\setlength{\unitlength}{1mm}
{\centering
\begin{tabular}{c|c|c|c} &\textbf{Before}&
\textbf{On Cool part} & \textbf{Beyond}\\\hline
\raisebox{7.5mm}[0pt]{$C_t$}&
\begin{picture}(32,20)
\put(0,0){\resizebox{32mm}{!}{\includegraphics{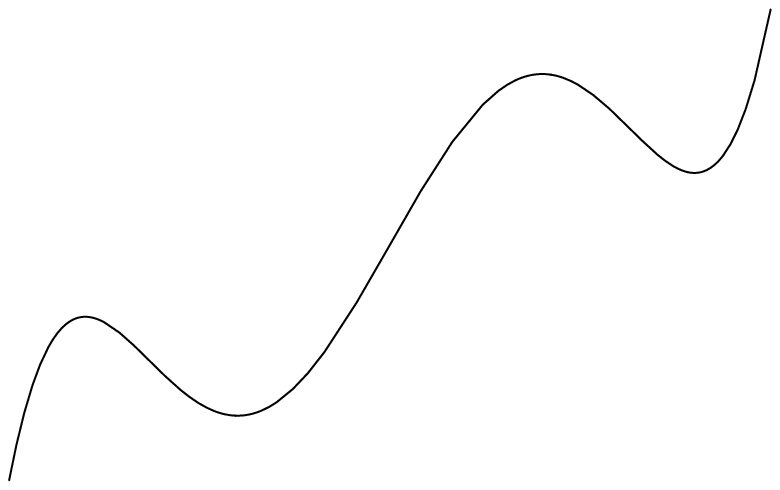}} }
\put(8,6){(a)}
\put(26,9){(b)}
\end{picture}&
\resizebox{32mm}{!}{\includegraphics{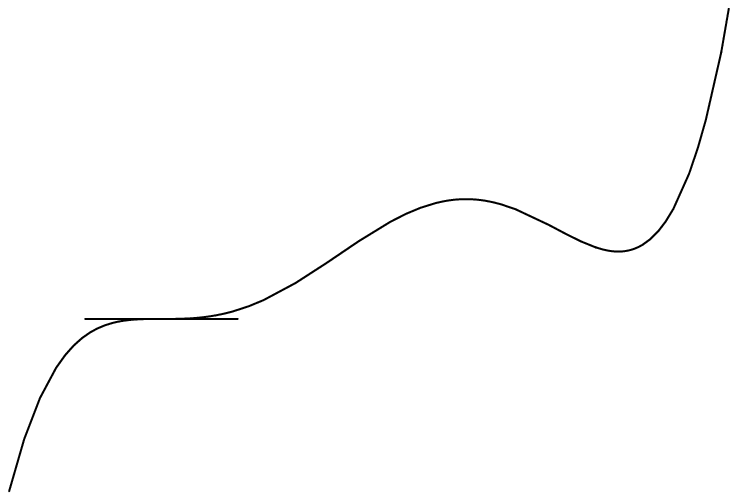}} &
\resizebox{32mm}{!}{\includegraphics{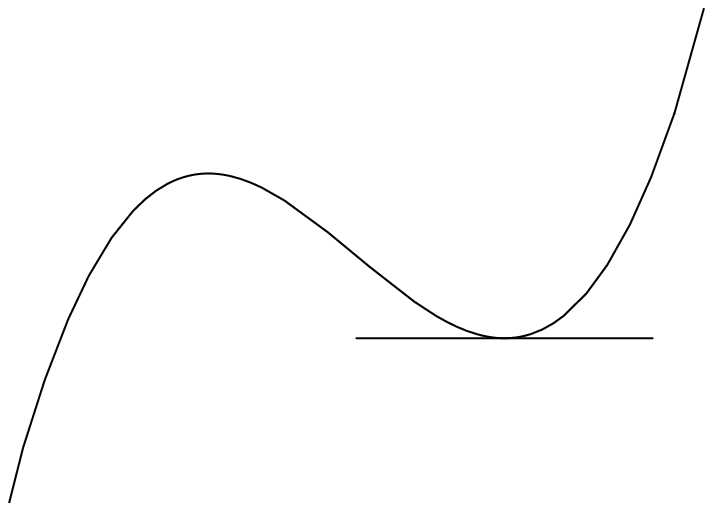}} \\\hline
\raisebox{7.5mm}[0pt]{$M_t$}&
\begin{picture}(32,20)
\put(0,0){\resizebox{32mm}{!}{\includegraphics{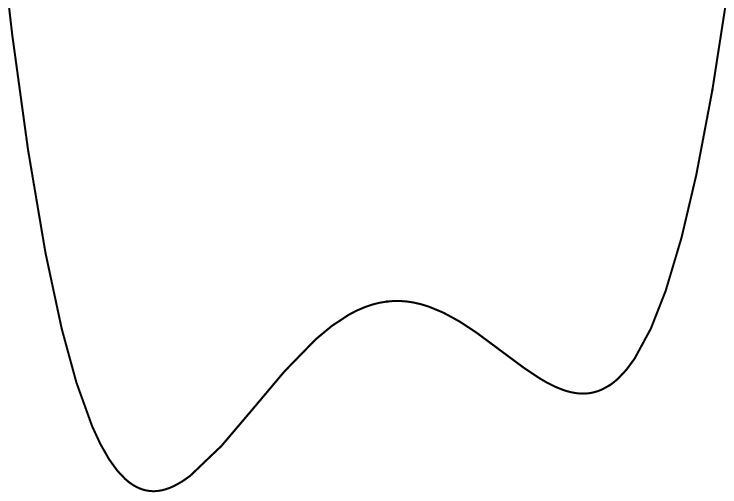}} }
\put(5.5,4){$x_0$}
\put(23,7){$\check{x}_0$}
\end{picture}&
\resizebox{32mm}{!}{\includegraphics{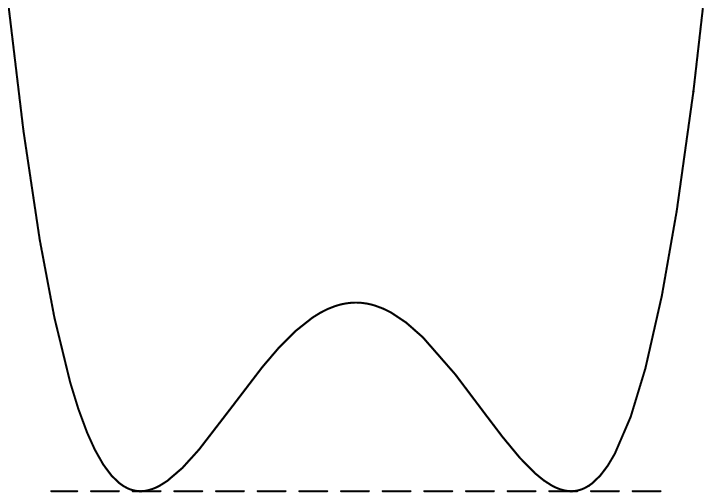}} &
\resizebox{32mm}{!}{\includegraphics{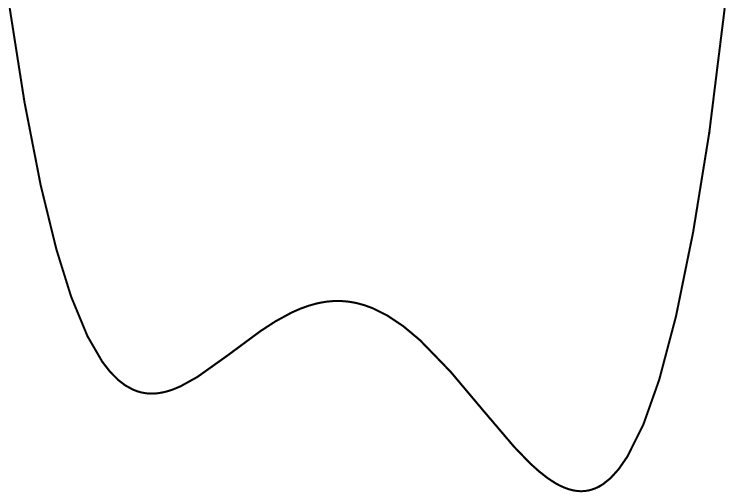}}
\end{tabular}
\caption{Graphs of the phase function as $x$ crosses $C_t$ and $M_t$.}
\label{fig1}}
\end{figure}

We illustrate this in one dimension by considering the integral,
 \begin{equation}\label{phase}I(x,t) = \int_{\mathbb{R}}G(x_0)\exp\left(\rmi
\frac{F(x_0,x,t)}{\mu^2}\right)\rmd x_0,\end{equation} where $G\in
C_0^{\infty}(\mathbb{R})$, $x\in\mathbb{R}^{d}$ and $\rmi=\sqrt{-1}$.
Consider the graph of the phase function, $F_{(x,t)}(x_0) =
F(x_0,x,t)$, as $x$ crosses the caustic and Maxwell set (see Figure \ref{fig1}).

 As we cross the caustic, the critical point at (a) becomes an inflexion which disappears causing $\tilde{x}_0(x,t)$ to jump from $(a)$ to $(b)$. This only
causes a jump in $v^{\mu}(x,t)$ when the point of inflexion is the global minimiser of $F$.
As we cross the Maxwell set, the critical points at $x_0$ and $\check{x}_0$ move so that $F_{(x,t)}(x_0)=F_{(x,t)}(\check{x}_0)$. If this
pair of critical points also minimise the phase function, then the
inviscid limit of the solution to the Burgers equation will jump.

\section{The reduced action function}

In this section we will find the phase function $F$ in equation (\ref{phase}).
We briefly summarise some results of Davies, Truman and Zhao \citep{MR1902482, MR2131267}.
As before, let the stochastic action be
defined as,
\[
A(x_0,p_0,t) =  \hf\int_0^t \dot{X}(s)^2\rmd s
 -\int_0^t\Bigg[ V(X(s))\rmd s +\epsilon k_s(X(s))\rmd
W_s\Bigg],\]
where
$X(s)=X(s,x_0,p_0)\in\mathbb{R}^d$ and for $ s\in[0,t]$ with $x_0,p_0\in\mathbb{R}^d$,
\[\rmd\dot{X}(s) = -\nabla V(X(s))\rmd s
-\epsilon \nabla k_s(X(s))\rmd W_s,\quad X(0)=x_0,
\quad\dot{X}(0)=p_0.\]
We assume $X(s)$ is unique and let $\mathcal{F}_s$ denote the sigma algebra generated by $X(u)$ up to time
$s$. It follows from Kunita \citep{MR776984}:

\begin{lemma} \label{i8}Assume $S_0,V\in
C^2$ and $k_t\in C^{2,0}$, $\nabla V,\nabla k_t$ Lipschitz with Hessians
$\nabla^2 V,\nabla^2 k_t$ and all second derivatives with respect to
space variables of $V$ and $k_t$ bounded. Then for $p_0$, possibly
$x_0$ dependent,
\[\frac{\partial A}{\partial
x_0^{\alpha}}(x_0,p_0,t) =\dot{X}(t)\cdot \frac{\partial
X(t)}{\partial x_0^{\alpha}}-\dot{X}_{\alpha}(0),\qquad\alpha
=1,2,\ldots,d.\]
\end{lemma}

The methods of
\citep{MR2073123}  guarantee that for small $t$ the map $p_0\mapsto
X(t,x_0,p_0)$ is onto for all $x_0$.
Therefore, we can define
$A(x_0,x,t)
:=A(x_0,p_0(x_0,x,t),t)$ where
$p_0=p_0(x_0,x,t)$ is the random minimiser (which we assume to be unique) of
$A(x_0,p_0,t)$ when $X(t,x_0,p_0)=x$.

Thus, the
stochastic action corresponding to the initial momentum $\nabla
S_0(x_0)$ is $\mathcal{A}(x_0,x,t) := A(x_0,x,t)+S_0(x_0).$

\begin{theorem}\label{sflow} If $\Phi_t$ is the stochastic flow map then,
\[\Phi_t(x_0) = x \quad\Leftrightarrow\quad\frac{\partial}{\partial x_0^{\alpha}}
\left[\mathcal{A}(x_0,x,t)\right]=0,\qquad \alpha=1,2,\ldots,d.\]
\end{theorem}

Using this we can create a one dimensional
\emph{reduced action function}. This is done by finding a series of
functions $x_0^{\alpha}(x_0^1,\ldots,x_0^{\alpha-1},x,t)$ for decreasing
$\alpha = d, d-1, \ldots,2$ by systematically locally solving the
equations,
\[\frac{\partial \mathcal{A}}{\partial x_0^{\alpha}}(x_0^1,\ldots,x_0^{\alpha},x_0^{\alpha+1}(\ldots),\ldots,x_0^d(\ldots),x,t)=0.\]
At each stage this eliminates one more coordinate from $x_0$ until
only $x_0^1$ remains. This gives \emph{local reducibility}
on the assumption that $\partial^2\mathcal{A}/(\partial
x_0^{\alpha})^2\neq 0$ for $\alpha =2,3,\ldots,d$ and also some mild
regularity conditions \citep{MR2083376}.
\begin{definition}
The reduced action function is the univariate function,
\[f_{(x,t)}(x_0^1)
:= \mathcal{A}(x_0^1,x_0^2(x_0^1,x,t),\ldots, x_0^d(x_0^1,x_0^2(\cdot),\ldots, x_0^{d-1}(\cdot),x,t),x,t).\]
\end{definition}

The Hamilton-Jacobi level surface $H_t^c$ is
found by eliminating $x_0$ between,
\[\mathcal{A}(x_0,x,t) =c, \qquad
\nabla_{x_0}\mathcal{A}(x_0,x,t) =0.\] Alternatively, if we eliminate $x$ to give
an expression in $x_0$, we have the pre-level surface
$\Phi_t^{-1}H_t^c$. Similarly the caustic $C_t$ (and pre-caustic
$\Phi_t^{-1}C_t$) are obtained by eliminating $x_0$ (or $x$)
between,
\[\det\left(
\frac{\partial^2\mathcal{A}} {\partial x_0^{\alpha}\partial
x_0^{\beta}}(x_0,x,t) \right)_{\alpha,\beta=1,2,\ldots, d} =0, \quad
\nabla_{x_0}\mathcal{A}(x_0,x,t) =0.\]
The Maxwell set $M_t$ (and pre-Maxwell set $\Phi_t^{-1} M_t$) are obtained by eliminating $x_0$ and $\check{x}_0$ (or $x$ and $\check{x}_0$) between the four equations,
$$\nabla_{x_0}\mathcal{A}(x_0,x,t) =0,\quad \nabla_{x_0}\mathcal{A}(\check{x}_0,x,t) =0, \quad \mathcal{A}(x_0,x,t)= \mathcal{A}(\check{x}_0,x,t)  = c.
$$

 The pre-images are calculated
algebraically and in the case of the pre-level surfaces are not necessarily the topological inverse
images. This can be done in the free case or when the relevant functions are polynomials in all variables which is an implicit assumption in what follows.

For polynomial $\mathcal{A}$, the eliminations involved with the Hamilton-Jacobi level surfaces and caustics are fairly simple to complete using the reduced action function with resultants and discriminants which can be calculated via Sylvester determinants \citep{Waerden}. The Maxwell set is more complicated to find as eliminating pre-images leads to a surface involving both real and complex pre-images termed the ``Maxwell-Klein set" \citep{MR2167533}. It is easier to find the pre-Maxwell set and then use the flow map to parameterise the Maxwell set. Parameterising in this manner allows one to restrict the pre-image of the Maxwell set to have only real values. In the polynomial case we have the following lemma,
\begin{lemma}
Let $D^x$ denote the polynomial discriminant taken with respect to $x$. The set of all singularities is,
\[D^c(D^{\lambda_1}(f_{(x,t)}(\lambda_1)-c))=0,\]
which factorises as,
\[k\times B_t(x)^2\times C_t(x)^3 = 0,\]
where $B_t=0$ is the equation of the Maxwell-Klein set, $C_t=0$ is the equation of the caustic and $k$ is some non-zero constant.

The pre-Maxwell set is given by,
\[D^{\lambda_1}\left(\frac{f_{(\Phi_t(x_0),t)}(x_0^1)-f_{(\Phi_t(x_0),t)}(\lambda_1)}
    {(x_0^1-\lambda_1)^2}\right)=0.\]
\end{lemma}

The reduced action function can also be used to identify the cool (singular) parts of the Maxwell set and caustic \citep{Progress}.

\section{Geometric Results}
The results in this section are taken from \citep{MR1902482, MR2167533, MR2314842}. Assume that $A(x_0,x,t)$ is $C^4$ in space variables with
$\det\left( \frac{\partial^2\mathcal{A}}{\partial
x_0^{\alpha}\partial x^{\beta}}\right)\neq 0.$
\begin{lemma}\label{s11}
Let $\Phi_t$ denote the stochastic flow map and
$\Phi_t^{-1}\Gamma_t$ and $\Gamma_t$ be some surfaces where if
$x_0\in\Phi_t^{-1}\Gamma_t$ then $x=\Phi_t (x_0)\in\Gamma_t$. Then,
$\Phi_t$ is a differentiable map from $\Phi_t^{-1}\Gamma_t$ to
$\Gamma_t$ with Frechet derivative,
\[(D\Phi_t)(x_0) = \left(-\frac{\partial^2\mathcal{A}}{\partial x\partial x_0}(x_0,x,t)\right)^{-1}
\left(\frac{\partial^2\mathcal{A}}{(\partial
x_0)^2}(x_0,x,t)\right).
\]
\end{lemma}
Let $n_{\mathrm{H}}(x_0)$, $n_{\mathrm{C}}(x_0)$ and $n_{\mathrm{M}}(x_0)$ denote the normal at $x_0$ to the pre-level surface, pre-caustic and pre-Maxwell set respectively. Using Lemma \ref{s11} we can show that:

\begin{theorem}\label{i11a}The normal to the pre-level
surface is, to within a scalar multiplier, given by,
\[n_{\mathrm{H}}(x_0) =
-\left(\frac{\partial^2\mathcal{A}}{(\partial x_0)^2}\right)
\left(\frac{\partial^2\mathcal{A}}{\partial x_0\partial x}\right)^{-1}
\dot{X}\left(t,x_0,\nabla S_0(x_0)\right).\]
\end{theorem}\begin{theorem}\label{m22}
Assume that a point $x$ on the Maxwell set corresponds to exactly
two pre-images on the pre-Maxwell set, $x_0$ and $\check{x}_0$. Then
the normal to the pre-Maxwell set at $x_0$ is, to within a scalar
multiplier, given by,
\begin{eqnarray*}
n_{\mathrm{M}}(x_0) &  = &  - \left(\frac{\partial^2\mathcal{A}}{(\partial
x_0)^2}(x_0,x,t)\right)
\left(\frac{\partial^2\mathcal{A}}{\partial x_0\partial x}(x_0,x,t)\right)^{-1}\cdot\\
& &\qquad \left(\dot{X}(t,x_0,\nabla
S_0(x_0))-\dot{X}(t,\check{x}_0,\nabla
S_0(\check{x}_0))\right)_.\end{eqnarray*}
\end{theorem}

We now consider the two dimensional case.
\begin{definition}
Let $x=x(\gamma)=(x_1,x_2)(\gamma)$ denote a curve where $\gamma$ is some intrinsic parameter (e.g. arc length)
with $\gamma\in (\gamma_0 - \delta, \gamma_0+\delta)$ for $\gamma_0\in\mathbb{R}$ and $\delta>0$. Then the curve is said to have a generalised cusp when $\gamma= \gamma_0$ if,
 \[\frac{\di x}{\di \gamma}( \gamma_0) =\left(\frac{ \di x_1}{\di \gamma}(\gamma_0),\frac{\di x_2}{\di \gamma}(\gamma_0)\right)=0.\]
\end{definition}
It then follows from Theorems \ref{i11a} and \ref{m22} that:
\begin{theorem}\label{h24}
Assume that in two dimensions at $x_0\in\Phi_t^{-1}H^c_t$ the normal
$n_{\mathrm{H}}(x_0)\neq 0$ so that the pre-level surface does not have a
generalised cusp at $x_0$. Then, the level surface can only have a
cusp at $\Phi_t(x_0)$ if $\Phi_t(x_0)\in C_t$. Moreover, if,
\[x=\Phi_t(x_0)\in\Phi_t\left\{ \Phi_t^{-1}C_t\cap
\Phi_t^{-1}H^c_t\right\},\] the level surface will have a generalised
cusp at $x$.
\end{theorem}
\begin{theorem}\label{m24}
Assume that in two dimensions at $x_0\in\Phi_t^{-1}M_t$ the normal
$n_{\mathrm{M}}(x_0)\neq 0$ so that the pre-Maxwell set does not have a
generalised cusp at $x_0$. Then, the Maxwell set can only have a
cusp at $\Phi_t(x_0)$ if $\Phi_t(x_0)\in C_t$. Moreover, if,
\[x=\Phi_t(x_0)\in\Phi_t\left\{ \Phi_t^{-1}C_t\cap
\Phi_t^{-1}M_t\right\},\] the Maxwell set will have a generalised
cusp at $x$.
\end{theorem}

These results lead to a range of conclusions relating the geometry of these curves. In particular, they allow us to characterise when swallowtails will form (a swallowtail perestroika). The appearance of a swallowtail is related to the existence of points with complex pre-images which are discussed in detail in \citep{MR2167533}.

\begin{corollary}\label{s14}
    Assume that at $x_0\in\Phi_t^{-1}H_t^c\cap\Phi_t^{-1}C_t$, $n_{\mathrm{H}}(x_0)\neq 0$ and  $n_{\mathrm{C}}(x_0)\neq 0$.
    Then at $\Phi_t(x_0)$ there is a cusp on the caustic if and only if $\Phi_t^{-1}H_t^c$ touches $\Phi_t^{-1}C_t$ at $x_0$. Moreover, it follows that, $x_0\in \Phi_t^{-1} M_t$ and that
    $\Phi_t^{-1}H_t^c$ touches $\Phi_t^{-1}M_t$ at $x_0$. Also, at $\Phi_t(x_0)$, $M_t$ will have a generalised cusp parallel to the cusp on $C_t$.
\end{corollary}
\begin{corollary}\label{m25}
 Assume that at $x_0\in\Phi_t^{-1}M_t\cap\Phi_t^{-1}C_t$, $n_{\mathrm{M}}(x_0)\neq 0$ and  $n_{\mathrm{C}}(x_0)\neq 0$.
 Then, there is a cusp on the Maxwell set where it intersects the caustic at
$x=\Phi_t(x_0)$ and the pre-Maxwell set touches a pre-level surface
$\Phi_t^{-1} H_t^c$ at $x_0$. Moreover, if the cusp on the Maxwell
set intersects the caustic at a regular point of the caustic, then
there will be a cusp on the pre-Maxwell set which also meets the
same pre-level surface $\Phi_t^{-1} H_t^c$ at another point
$\check{x}_0$.\end{corollary}

\begin{corollary}\label{s14a}
    Assume that at $x_0\in\Phi_t^{-1}H_t^c\cap\Phi_t^{-1}C_t$, $n_{\mathrm{H}}(x_0)\neq 0$ and  $n_{\mathrm{C}}(x_0)\neq 0$.
    Then at $\Phi_t(x_0)$ there is a point of swallowtail perestroika on the level surface
    $H_t^c$ if and only if  there is a generalised cusp on the caustic $C_t$ at $\Phi_t(x_0)$.
\end{corollary}

The results in this section have natural extensions to three dimensions where the cusps are replaced by curves of cusps.
We now give some two dimensional examples.

\begin{eg}[The generic cusp]
 We consider a two dimensional deterministic free example ($V \equiv 0$, $\epsilon = 0$). In general for such a system the flow map is given by,
\[\Phi_t(x_0) = x_0 + t\nabla S_0(x_0),\]
with derivative map $D \Phi_t (x_0) = (I + t\nabla^2S_0(x_0)).$
The  pre-level surface is then given by the eikonal equation,
\[\frac{t}{2}|\nabla S_0(x_0)|^2 +S_0(x_0) = c,\]
 where the key identity is,
\[\nabla_{x_0}\left\{\frac{t}{2}|\nabla S_0(x_0)|^2 +S_0(x_0) \right\} = (I + t\nabla^2S_0(x_0))\nabla S_0(x_0).\]
The generic cusp initial condition, $S_0(x_0,y_0) =  x_0^2y_0,$ gives a simple cusped caustic (see Figure \ref{gencusp}).
\begin{figure}[t]
\setlength{\unitlength}{8mm}
  \begin{tabular}{cc}
 \begin{picture}(7,7)
\put(0,0){\resizebox{56mm}{!}{\includegraphics{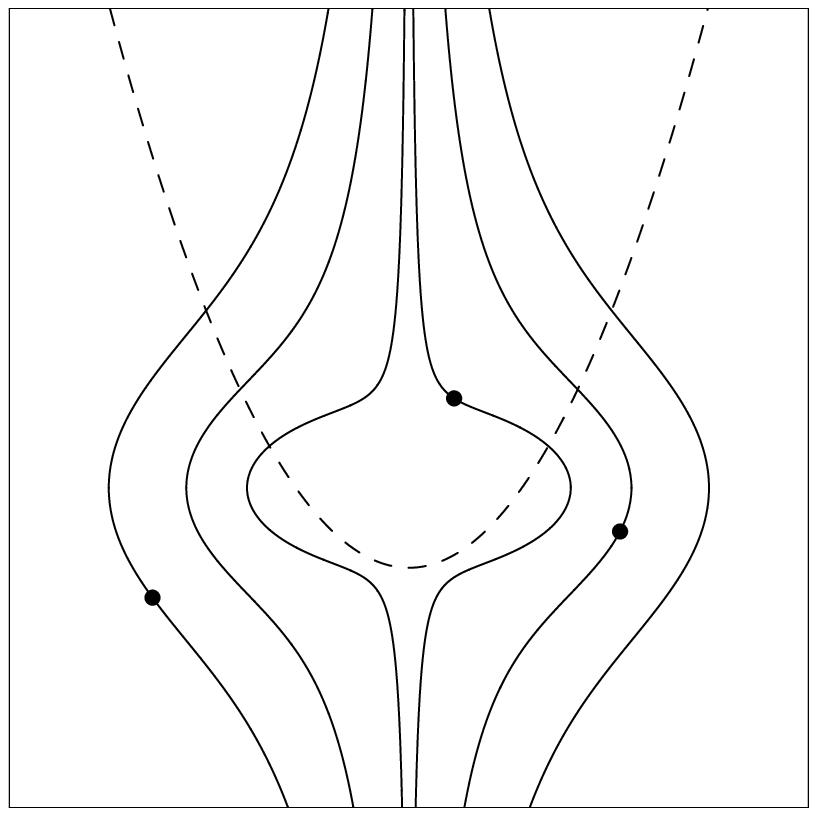}}}
 \end{picture}&
\begin{picture}(7,7)
 \put(0,0){\resizebox{56mm}{!}{\includegraphics{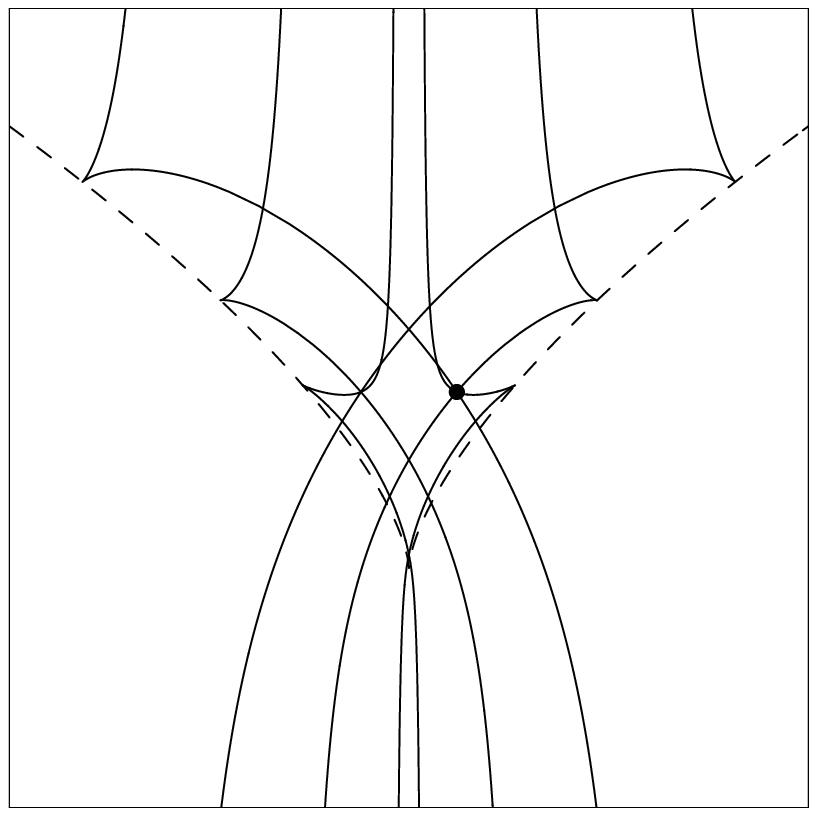}}}
 \end{picture}\\
 $\Phi_t^{-1}C_t$ and $\Phi_t^{-1}H^c_t$& $C_t$ and $H_t^c$
\end{tabular}
\caption{The generic cusp caustic (dashed) with three level surfaces (solid line)}
\label{gencusp}
\end{figure}
 \end{eg}

\begin{eg}[The polynomial swallowtail]
Let $V(x,y)=0$, $k_t(x,y)=x$ and 
 $S_0(x_0,y_0)=x_0^5+x_0^2y_0.$
The noisy potential does not affect either the pre-caustic or pre-Maxwell set. Consequently at time $t$ the noise will have shifted the deterministic caustic and Maxwell set by $-\epsilon \int_0^t W(u) \di u$ in the $x$ direction. This point will be returned to in Section 6.

\begin{figure}[p]
\setlength{\unitlength}{8mm}
\centering \begin{tabular}{cc}
 \begin{picture}(7,7)
 \put(2.2,4.6){1}
  \put(0.4,5.2){2}
 \put(1.7,5.8){3}
  \put(4.8,5.3){4}
  \put(6.4,0.2){5}
 \put(4.8,3){6}
\put(0,0){\resizebox{56mm}{!}{\includegraphics{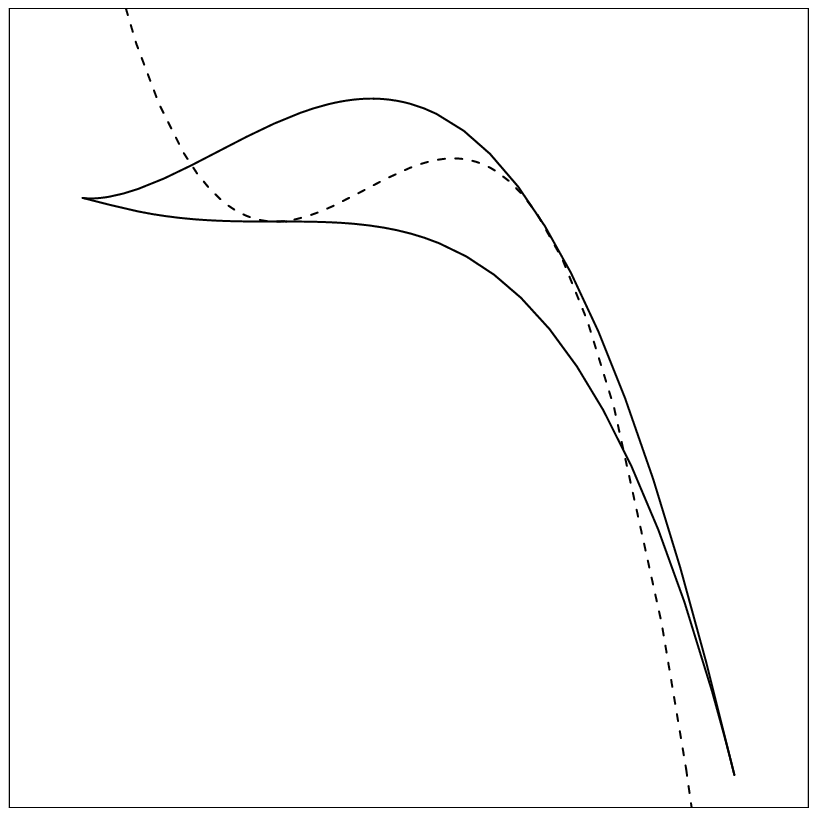}}}
\end{picture}
&
 \begin{picture}(7,7)
 \put(2.6,0.2){1}
  \put(3.1,4.5){2}
 \put(2,2.6){3}
  \put(6.7,6.3){4}
  \put(2.1,2.2){5}
 \put(3.4,4.6){6}
\put(0,0){\resizebox{56mm}{!}{\includegraphics{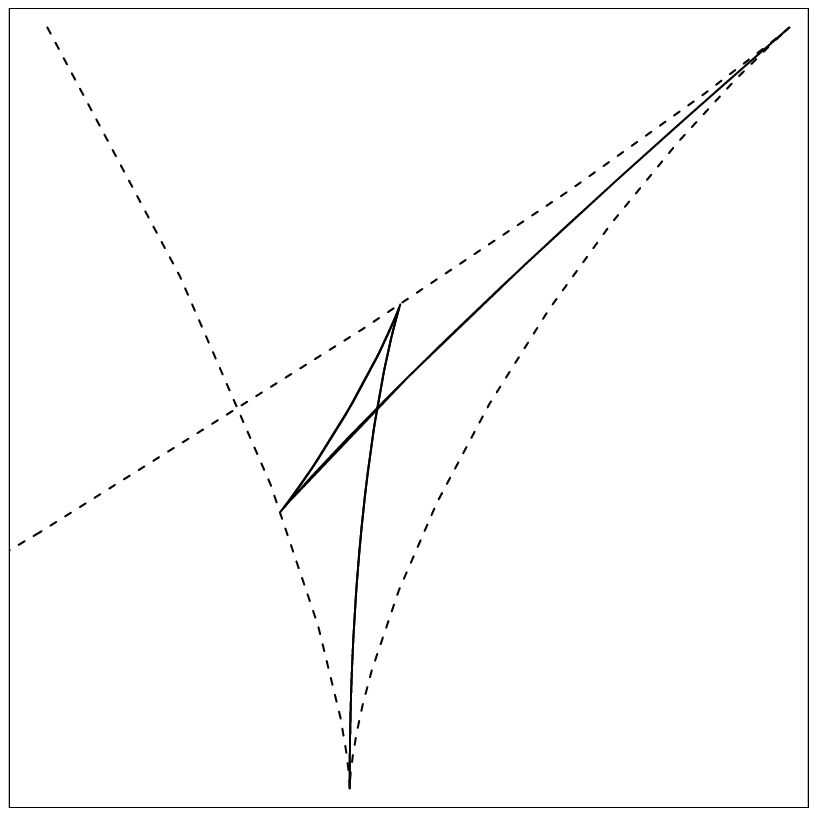}}}
\end{picture}
\\
\emph{ Pre-curves}   & \emph{ Curves } \\
\end{tabular}
\caption{The polynomial swallowtail  caustic (dashed) and Maxwell set (solid line).} \label{mf4}
\end{figure}
From Theorem \ref{m24}, the cusps on the Maxwell set correspond to the intersections of the pre-curves (points 3 and 6 on Figure \ref{mf4}). But from Corollary \ref{m25}, the cusps on the Maxwell set also correspond to the cusps on the pre-Maxwell set (points 2 and 5 on Figure \ref{mf4} and also Figure \ref{mf4a}). The Maxwell set terminates when it reaches the cusps on the caustic. These points satisfy the condition for a generalised cusp but, instead of appearing cusped, the curve stops and maps back exactly onto itself. At such points the pre-surfaces all touch (Figure \ref{mf4a}).

These two different forms of cusps correspond to very different geometric behaviours of the level surfaces. From the definition of a Maxwell set it is clear that any point on $M_t$ is a point of self intersection of some level surface. Where the Maxwell set stops or cusps corresponds to the disappearance of a point of self-intersection on a level surface. There are two distinct ways in which this can happen. Firstly, the level surface will have a point of swallowtail perestroika when it meets a cusp on the caustic. At such a point only one point of self-intersection will disappear, and so there will be only one path of the Maxwell set which will terminate at that point. However, when we approach the caustic at a regular point, the level surface  must have a cusp but not a  swallowtail perestoika. This corresponds to the collapse of two points of self intersection and
so two paths of the Maxwell set must approach the point and produce the cusp
(see Figure \ref{mf5}).

\begin{figure}[p]
\centering \begin{tabular}{cc}
\resizebox{56mm}{!}{\includegraphics{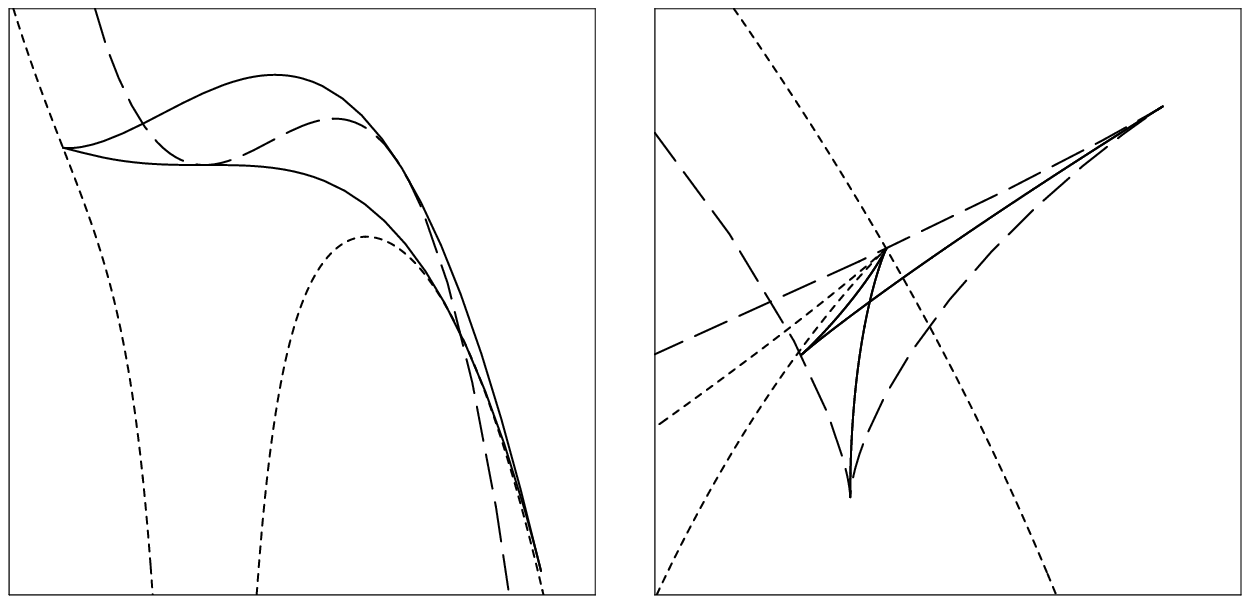}}&
\resizebox{56mm}{!}{\includegraphics{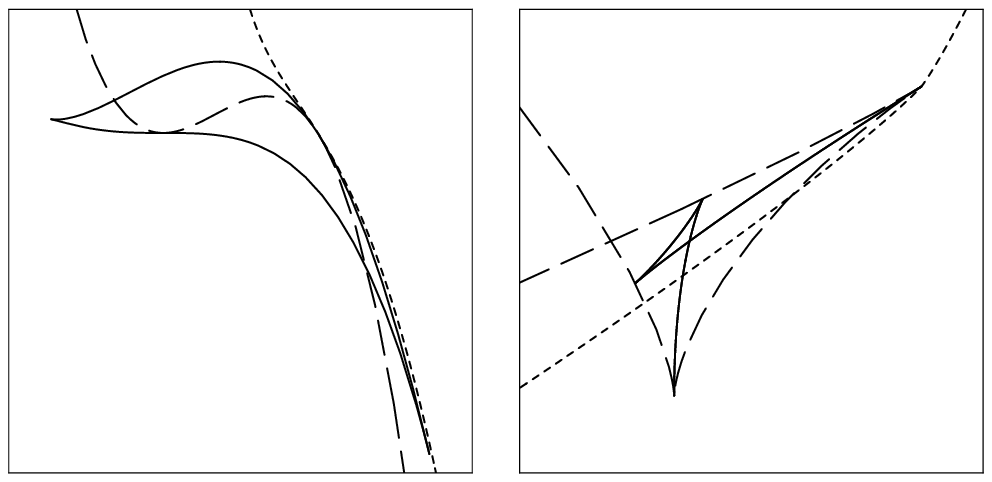}}\\
\emph{ Cusp on Maxwell set}   & \emph{ Cusp on caustic } \\
\end{tabular}
\caption{The caustic (long dash) and Maxwell set (solid line) with the level surfaces (short dash) through special points.} \label{mf4a}
\end{figure}

\begin{figure}[p]\centering
\begin{tabular}{cccc}
\resizebox{26mm}{!}{\includegraphics{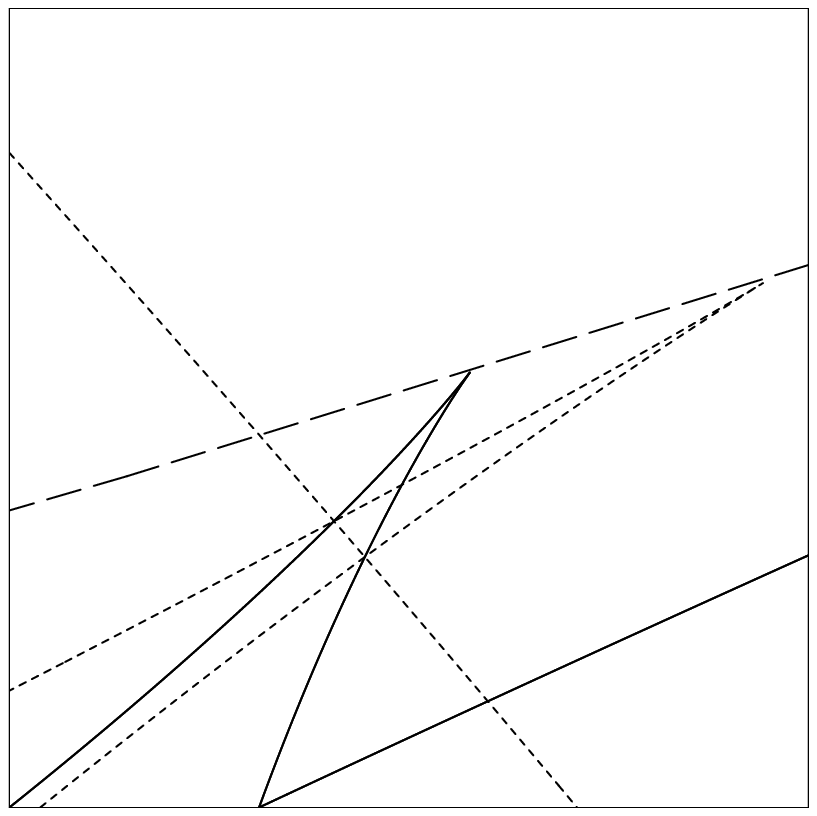}}&
\resizebox{26mm}{!}{\includegraphics{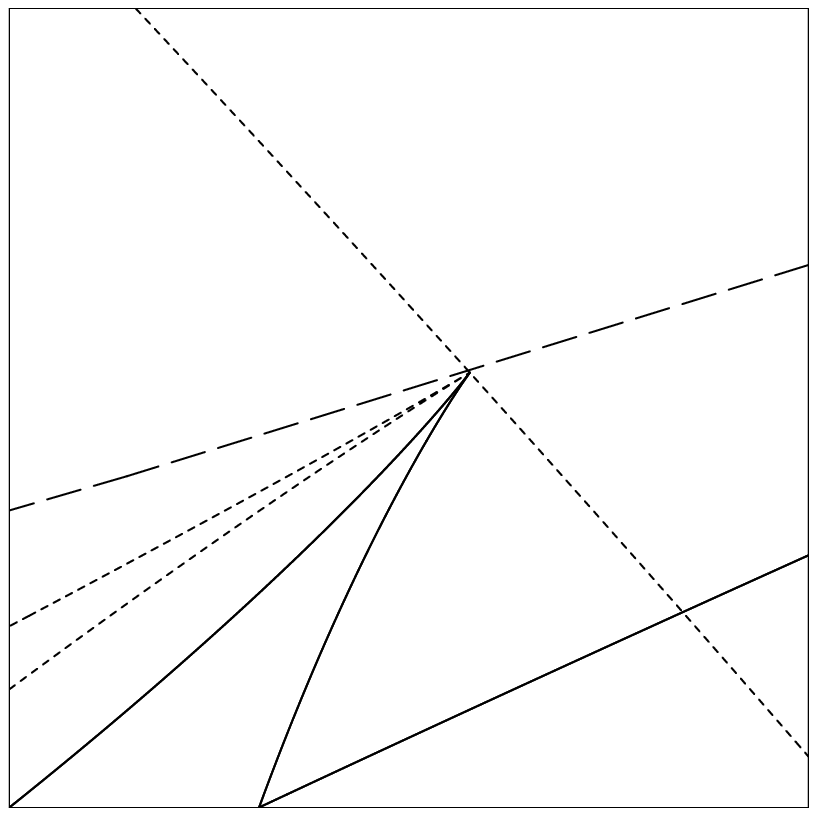}}&
\resizebox{26mm}{!}{\includegraphics{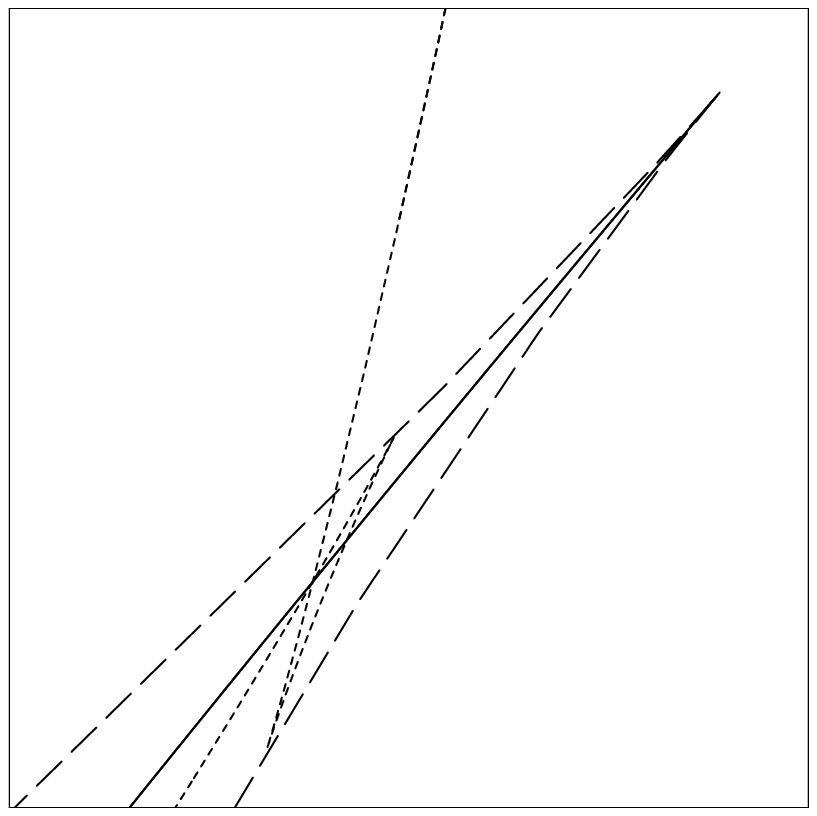}}&
\resizebox{26mm}{!}{\includegraphics{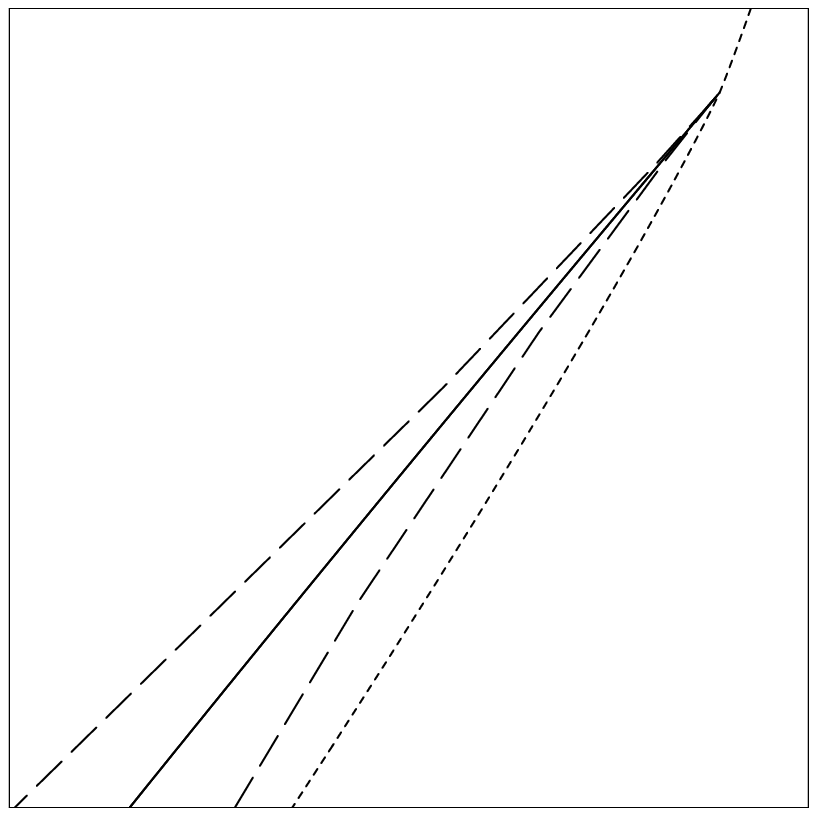}}
\\
\multicolumn{2}{c}{\emph{ Approaching the caustic } }  & \multicolumn{2}{c}{\emph{ Approaching a cusp on the caustic }}
\end{tabular}
\caption{The caustic (long dash) Maxwell set (solid line) and level surface (short dash).} \label{mf5}
\end{figure}
\end{eg}

\section{Recurrence of stochastic turbulence}

Following   \citep{MR2083376, mprf}, the geometric results of Section 5  can be used to characterise a sequence of turbulent times.

\begin{definition}\label{t1}
    Real turbulent times are defined to be times $t$ at which there
    exist real points where the
    pre-level surface $\Phi_t^{-1}H_t^c$ and pre-caustic $\Phi_t^{-1} C_t$ touch.
\end{definition}
Real turbulent times correspond to times at which there is a change
in the number of cusps or cusped curves on the level surface
$H_{t}^{c}$. In $d$-dimensions, assuming $\Phi_t$ is
globally reducible, let $f_{({x},t)}(x_{0}^1)$ denote the reduced
action function and $x_t(\lambda)$ denote the caustic parameterised using the pre-caustic and flow map.

\begin{theorem}\label{t2}
The real turbulent times $t$ are given by the zeros of the zeta
process $\zeta_t^c$ where,
\[\zeta_t^c:=f_{(x_{t}(\lambda),t)}(\lambda_1)-c,\]
$\lambda$ satisfies,
\begin{equation}\label{te1b}\frac{\partial}{\partial \lambda_{\alpha}}f_{(x_{t}(\lambda),t)}(\lambda_1)
=0\quad\mbox{for }\alpha=1,2,\ldots,d,\end{equation}
and $x_t(\lambda)$ is on the cool part of the caustic.
\end{theorem}
The term `real' is used in \citep{MR2167533} to distinguish this form of turbulence from `complex' turbulence where swallowtail perestroikas occur on the caustic. We shall not discuss the details of complex turbulence in this article.

We now consider the stochastic Burgers equation with white noise forcing in $d$-orthogonal directions,
\begin{equation}\label{te1a}\frac{\partial v^{\mu}}{\partial t}
+(v^{\mu}\cdot\nabla)v^{\mu} = \hf[\mu^{2}]\Delta v^{\mu}
-\epsilon {\dot{W}}(t),\end{equation}
where  $W(t)=\left(W_1(t),W_2(t),\ldots,W_d(t)\right)$ is a $d$-dimensional Wiener process.
\begin{proposition}\label{t3}
The stochastic action corresponding to the Burgers equation (\ref{te1a})
is,
\begin{eqnarray*}
\mathcal{A}({x_{0}},{x},t)
& = &
\frac{|x-x_0|^{2}}{2t}
+\frac{\epsilon}{t}(x-x_0)\cdot\int_{0}^{t}W(s)\di s
-\epsilon x\cdot W(t)
\\
&   & \quad
-\hf[\epsilon^{2}]\int_{0}^{t}|W(s)|^2\di s+\frac{\epsilon^{2}}{2t}\left|\int_{0}^{t}W(s)\di
u\right|^{2}
+S_{0}({x_{0}}).
\end{eqnarray*}
\end{proposition}

\begin{lemma}\label{t6}
 If $x_t^{\epsilon}(\lambda)$ denotes the random caustic for the stochastic Burgers equation (\ref{te1a}) and $x_t^0(\lambda)$ denotes the deterministic caustic (the $\epsilon=0$ case) then,
\[x_t^{\epsilon}(\lambda)=x_t^0(\lambda)-\epsilon\int_0^t W(u) \di u.\]
\end{lemma}
Using Proposition  \ref{t3} and Lemma \ref{t6}, we can find the zeta process explicitly.
\begin{theorem}\label{t7}
In $d$-dimensions, the zeta process for the stochastic Burgers equation (\ref{te1a}) is,
\[
\zeta_t^c  = f_{(x_t^0(\lambda),t)}^0(\lambda_1) -\epsilon
x_t^0(\lambda)\cdot W(t) +\epsilon^2 W(t)\cdot \int_0^t W(s)\di s
-\hf[\epsilon^2]\int_0^t |W(s)|^2 \di s-c,\]
where $f_{(x,t)}^0(\lambda_1)$ is the deterministic reduced action
function, $x_t^0(\lambda)$ is the deterministic caustic and
$\lambda$ must satisfy the stochastic equation,
\begin{equation}\label{te5a}\nabla_{\lambda}\left(f^0_{(x_t^0(\lambda),t)}(\lambda_1)-
\epsilon x_t^0(\lambda)\cdot W(t)\right)
=0.\end{equation}
\end{theorem}
Equation (\ref{te5a}) shows that the value of $\lambda$ used in the zeta process may be either deterministic or random.
In the two dimensional case this gives,
\begin{equation}
0 =
\left( \nabla_x f^0_{(x_t^0(\lambda),t)}(\lambda_1)-
\epsilon  W(t)\right)\cdot\frac{\di x_t^0}{\di\lambda}(\lambda),\label{te5b}
\end{equation}
which  has a deterministic solution for $\lambda$ corresponding to a
cusp on the deterministic caustic. 

Using the law of the iterated logarithm, it is a simple matter to
show formally that if there is a time $\tau$ such that
$\zeta_{\tau}^c=0$, then there will be infinitely many zeros of
$\zeta_t^c$ in some neighbourhood of $\tau$. This suggests that the set
of zeros of $\zeta_t^c$ a perfect set which can be rigorously proved in some generality \citep{Chris}.

The intermittence of turbulence will be demonstrated if we can show
that there is an unbounded increasing sequence of times at
which the zeta process is zero.
This can be done using an idea of David Williams and the Strassen form of the law of the iterated logarithm  \citep{MR2083376}.
\begin{theorem}\label{t12}
There exists an unbounded increasing sequence of times $t_n$ for which $Y_{t_n}=0$, almost surely, where,
\[
Y_t =   W(t)\cdot \int_0^t W(s)\di s
-\hf\int_0^t |W(s)|^2 \di s,
\]
and $W(t)$ is a $d$-dimensional Wiener process.
\end{theorem}

\begin{corollary}\label{t13}
Let $h(t) = (2 t \ln\ln(t))^{-\hf}$. If
$h(t)^2t^{-1}f^0_{(x_t^0(\lambda),t)}(\lambda_1)\rightarrow 0$ and $h(t)t^{-1}\sum\limits_{i=0}^d x_t^{0_i}(\lambda)\rightarrow 0,$
then the zeta process $\zeta_t^c$ is recurrent.
\end{corollary}

\section{A vortex line sheet on the Maxwell set}
We now summarise results of \citep{MR2314842}  which show that for compressible
flow and small viscosity, with appropriate initial conditions, there
is a vortex filament structure in the neighbourhood of the cool part
of the Maxwell set. This result is valid for both
deterministic and stochastic cases.

We first compute the Burgers
fluid velocity on the cool part of the Maxwell set where typically
for $x\in M_t$ the pre-images $x_0(x,t)$ and $\check{x}_0(x,t)$ are well behaved
functions.
Recall from equations (\ref{eqn3}) and (\ref{eqn2}) that for $t<t_c(\omega)$ and small $\mu$,
\[u^{\mu}(x,t) \sim \exp\left(-\frac{\mathcal{S}_t(x)}{\mu^2}\right)
\rho_t^{\hf}(x) (1+\Or(\mu^2)).\]
For $t>t_c(\omega)$, the analagous result for $v^{\mu}$ is:

\begin{lemma}\label{V1}
Let $x\in\mathrm{Cool}(M_t)$, so that
$x=\Phi_t(x_0)=\Phi_t(\check{x}_0)$ where $x_0\neq \check{x}_0$ and
$\mathcal{A}(x_0,x,t) = \mathcal{A}(\check{x}_0,x,t)$. Then,
\[v^{\mu}(x,t)\sim \frac{\rho_t^{\hf}(x)\nabla \mathcal{S}_t(x)
+\check{\rho}_t^{\hf}(x)\nabla
\check{\mathcal{S}}_t(x)}{\rho_t^{\hf}(x)+\check{\rho}_t^{\hf}(x)}+\mathrm{O}(\mu^2),\]
where,
\[\rho_t^{\hf}(x) = T_0(x_0(x,t))\left|\frac{\partial x_0}{\partial x}(x,t)\right|^{\hf}, \quad
\check{\rho}_t^{\hf}(x) = T_0(\check{x}_0(x,t))\left|\frac{\partial
\check{x}_0}{\partial x}(x,t)\right|^{\hf},\] and
$\mathcal{S}_t(x)=\mathcal{A}(x_0(x,t),x,t),$ $\check{\mathcal{S}}_t(x)
= \mathcal{A}(\check{x}_0(x,t),x,t).$
\end{lemma}

We denote by $v^0(x,t)$ the leading term for the behaviour of
$v^{\mu}(x,t)$ on  the Maxwell set and choose orthogonal curvilinear coordinates on $M_t$ denoted by
$(\xi_1,\xi_2)$. Let the unit normal in a coordinate patch on $M_t$ be denoted $n$.

\begin{theorem}\label{v1} If $v^0(x,t)$  is the leading behaviour of
$v^{\mu}(x,t)$ for $x\in M_t$, then,
\[v^0(x,t) =\hf\left\{
    \nabla(\mathcal{S}_t(x)+\check{\mathcal{S}}_t(x)) +
        \frac{\check{\rho}_t^{\hf}(x)-{\rho}_t^{\hf}(x)}{\check{\rho}_t^{\hf}(x)+{\rho}_t^{\hf}(x)}
            \left(
                \frac{\partial \check{\mathcal{S}}_t}{\partial n}(x)-\frac{\partial {\mathcal{S}}_t}{\partial n}(x)\right)n\right\}_,\]
                where $\frac{\partial}{\partial n}$ denotes the normal derivative $(n\cdot\nabla)$ on $M_t$.
\end{theorem}

\begin{definition}
The inviscid limit of the vorticity $\omega^0$ is defined to be,
     \[\omega^0(x,t) := \nabla\wedge v^0(x,t),\]
     where $v^0(x,t)$ denotes the leading behaviour of $v^{\mu}(x,t)$ as $\mu\rightarrow 0$.
\end{definition}
Since the first term in $v^{0}$ is $C^1$ on $M_t$, we obtain:
\begin{corollary}\label{kelvin}
For $x\in\mathrm{Cool}(M_t)$,
\[\omega^0 = \nabla\wedge\left\{
    \frac{\check{\rho}_t^{\hf}(x)-{\rho}_t^{\hf}(x)}{\check{\rho}_t^{\hf}(x)+{\rho}_t^{\hf}(x)}
            \left(
                \frac{\partial \check{\mathcal{S}}_t}{\partial n}(x)-\frac{\partial {\mathcal{S}}_t}{\partial n}(x)\right)n\right\}
                    \in T_xM_t. \]
\end{corollary}

From Corollary \ref{kelvin} it follows that $\omega^{0}\ne 0$. Hence, we expect that for small viscosity even though inititially there was zero vorticity, once a Maxwell set appears, the flow is no longer irrotational.

The above result shows that for small viscosity a
vortex filament structure will appear in a neighbourhood of the cool
part of the Maxwell set. The limit of this vortex filament structure
is a sheet of vortex lines on the cool part of the Maxwell set. We
now give the equation of these limiting vortex lines on the
Maxwell set in terms of orthogonal coordinates $(\xi_1,\xi_2)$.

\begin{theorem}
The limiting vortex lines in a coordinate patch of $M_t$ have
equations,

\[h_1(\xi)h_2(\xi)(\check{\rho}_t^{\hf}(x)-\rho_t^{\hf}(x))\left(
\frac{\partial}{\partial n}\left( f_{(x,t)}(\check{x}_0^1(x,t))
-f_{(x,t)}({x}_0^1(x,t)) \right)\right)\]
\[\qquad\qquad\qquad\qquad\qquad=c(\check{\rho}^{\hf}_t(x)+\rho^{\hf}_t(x)),
\]
where $c$ is a real constant and $\frac{\partial}{\partial n}$
denotes the normal derivative $(n\cdot\nabla)$ on $M_t$.
\end{theorem}

This confirms that for our initial conditions, for small
viscosity, and for compressible flow, vortex filaments will inevitably appear in a neighbourhood of the
cool part of the Maxwell set.  Given the rotational effects at work in the universe perhaps this suggests that we should consider the Burgers equation with vorticity from the
outset. Kinematical considerations and Galilean invariance suggest
that the appropriate equation is a Burgers equation with a vector
potential. We hope to discuss this in a future paper \citep{Vortex}.

\section{The adhesion model}
The adhesion model for the formation of the early universe is a refinement of the Zeldovich approximation \citep{MR657014}.
In the original adhesion model  there is no noise and a variational principle is assumed which forces the mass to move perpendicular to the cool part of the Maxwell set (often referred to as the shock) with the same velocity as the Maxwell set itself. This clearly results in mass adhering to the cool Maxwell set leading to an accumulation of mass at certain points \citep{Bec, Bog}.  A simple calculation gives the velocity of the Maxwell set as,
\[\hf\left\{
    \nabla(\mathcal{S}_t(x)+\check{\mathcal{S}}_t(x)) \right\}_,\]
and a comparison with Theorem \ref{v1} and Corollary \ref{kelvin}, reveals that this adhesion is sufficient to precisely  destroy the vorticity near the Maxwell set.
When adhesion occurs, an interesting question to consider is how rapidly mass will accrete on the Maxwell set as $\mu\rightarrow 0$ \citep{Bec}. This would relate directly to the mass involved in the formation of galaxies on the shock. Here we prove an inequality 
for the magnitude of the accumulated mass.

For this analysis we need to first mollify the Nelson diffusion process introduced in Section 2, equation (\ref{ItoX}),
\[ \di X_u^{\mu} = \nabla \mathcal{S}_u (X^{\mu}_u) \di u +\mu \di B(u),\quad X_t^{\mu} = x.\]
This mollification will remove the discontinuities in the drift caused by the Maxwell set and caustic.
 
Let  $0<t<T$, for some fixed $T$,  and let $u\in(0,t)$.  Assume that we can mollify the minimising Hamilton-Jacobi function $\mathcal{S}_u$ such that,
   $\mathcal{S}_u = \mathcal{S}_u^{\mathrm{moll}}$ off some thin open set $\tau_u(\mu)$ surrounding the cool Maxwell set $M_u^{\mathrm{Cool}}$ and the cool caustic $C^{\mathrm{Cool}}_u$ where  the Lebesgue measure $|\tau_u(\mu)| = \mathrm{O}(\mu)$, and that $\nabla \mathcal{S}_u^{\mathrm{moll}}$ is uniformly Lipschitz in space and bounded for $u\in(0,t)$.

We can then define a mollified potential, $V_u^{\mathrm{moll}}(x)$, corresponding to this new system  such that,
 \[\frac{\partial \mathcal{S}_u^{\mathrm{moll}}}{\partial u}+\hf |\nabla \mathcal{S}_u^{\mathrm{moll}}|^2+ V^{\mathrm{moll}}_u = 0,\]
so that as $\mu\rightarrow 0$, $(V_u^{\mathrm{moll}}-V)$ is a surface potential.

 This system then has a corresponding Burgers equation,
 \[\frac{D v^{\mathrm{moll}}}{Dt} = \hf[\mu^2] \Delta v^{\mathrm{moll}} - \nabla V^{\mathrm{moll}},\]
 where this is the physically important velocity field in the limit $\mu\rightarrow 0$.

 Let $\mu = \mu_n$ where $\mu_n$ is a sequence of real values such that $\mu_n\rightarrow 0$ as $n\rightarrow \infty$.
 Let the diffusion associated with $\mu  = \mu_n$ after mollification be denoted by $X^n$ and let, 
 \[A_t^n = \left\{ \omega: X^n(\omega) \mbox{ avoids  } C^{\mathrm{Cool}}_u \cup M^{\mathrm{Cool}}_u  \mbox{ at all times }u\in(0,t)\right\}. \]
It is important to note that, formally at least,
\[\mathbb{P}(A_t^n) = \lim\limits_{\lambda\rightarrow\infty} \mathbb{E}_x\left[ \exp \left(-\lambda \int_0^t\chi_{C_u^{\mathrm{Cool}}(\lambda^{-1})}(X^n_{\mathrm{trev}}(u)) \di u\right)\right.\]
\[\qquad\qquad\left. \times\exp \left(-\lambda \int_0^t\chi_{M_u^{\mathrm{Cool}}(\lambda^{-1})}(X^n_{\mathrm{trev}}(u)) \di u\right)\right],\]
where,
 \begin{eqnarray*}
 C_s^{\mathrm{Cool}}(\lambda^{-1})&:= &\left\{y: d(y,C_s^{\mathrm{Cool}}) <\lambda^{-1}\right\},\\
  M_s^{\mathrm{Cool}}(\lambda^{-1})&:= &\left\{y: d(y,M_s^{\mathrm{Cool}}) <\lambda^{-1}\right\}.
\end{eqnarray*}
Here trev denotes time reversal so that,
\[ X^n_{\mathrm{trev}} (s) := X^n(t-s),\qquad X^n_{\mathrm{trev}}(0)=x.\]
We need to assume that for sufficiently large $n$, $\mathbb{P}(A_t^n) >\delta>0.$

Now assume that at any point $x$ on the cool Maxwell set, $M_t^{\mathrm{Cool}}$ divides space into $m(x)$ parts.
Therefore, $x$ has $m(x)$ minimising pre-images  $x_0(i)(x,t)$, ($i = 1, 2,\ldots,m$)
and for each pre-image there is a corresponding $S = S(i)$ given by,
\[\nabla S_u (i)(x) =\left. \dot{X}(x_0,\nabla S_0(x_0),u)\right|_{x_0 = x_0(i)(x,t).}\]

We assume that all $\nabla S(i)$ are bounded and uniformly Lipschitz in space away from the caustic.
Note that at a regular point of the Maxwell set $m=2$ which is the case in which we are interested. In this case $x_0(i)(x,t)$ for $i=1,2$ correspond to matter arriving on different sides of the cool part of the Maxwell set.
Arguing as in \citep{MR1652127}, the Borel-Cantelli and Gronwall's lemmas give:

\begin{lemma}
 For $u\in(0,t)$ and $n\in\mathbb{N}$ ($i=1,2$) let,
 \[\di X_i^n(u) = \nabla S^n_u(i)(X_i^n(u))\di u+\mu_n\di B(u), \]
 and
 \[\di X_i^0(u) = \nabla S_u (i)(X_i^0(u)) \di u,\]
 with $X_i^n(t) = X_i^0(t) = x$. Then, if $\sum \mu^2_n<\infty,$ for $0<t<T$,
 \[\mathbb{P}\left[\sup\limits_{0<u<t}\left|X_i^n(u) - X_i^0(u)\right|\rightarrow 0 \mbox{ as } n\rightarrow \infty \bigg|A_t^n\right] = 1.\]
\end{lemma}

We can now give our results on the distribution of mass.  We parameterise the cool Maxwell set using the pre-Maxwell set and flow map,
 \[x_0 = (x_0^1,x_0^2,\ldots,x_0^{d-1},x_0^d(x_0^1,\ldots,x_0^{d-1},t))\in\Phi_t^{-1}M^{\mathrm{Cool}}_t.\]
 Define $x_0 \in \left(\Phi_t^{-1} M_t^{\mathrm{Cool}}\right)^{\tilde{}} $
if  the classical path from $x_0$ to $M_t^{\mathrm{Cool}}$ avoids $C_u^{\mathrm{Cool}}$ and $M_u^{\mathrm{Cool}}$ for all $u<t$.
It now follows from our results in Section 2 that:

\begin{theorem}
Let the mass adhering to the Maxwell set in time interval $(0,T)$ be $m(0,T)$. Then if,
\[m_0(0,T) = \int_0^T \hf[\di t ]\int_{x_0\in(\Phi_t^{-1}M_t^{\mathrm{Cool}})^{\tilde{}}} T_0^2(x_0) \left|\frac{\partial x_0^d}{\partial t}\right| \di x_0^1\ldots \di x_0^{d-1},\]
we obtain,
\[m(0,T) \ge m_0(0,T).\]
\end{theorem}

 The volume of $M_t^{\mathrm{Cool}}$ is zero and $m(0,T)$ is O(1). Therefore, the mass density on $M_t^{\mathrm{Cool}}$ will be infinite in parts of $M_t^{\mathrm{Cool}}$ in the inviscid limit. As we stated earlier, the process of adhesion forces the velocity away from $M_t$ to be zero, and consequently the adhering particles carry no vorticity.
 
We now conjecture what happens when particles hit $C_t^{\mathrm{Cool}}$ first.
Taking limits, $m_0(0,T)$ is the contribution to the mass in the shock from paths with no kinks which would be caused by hitting the caustic before time $t$.

It would therefore be reasonable to calculate $m_1(0,T)$,  the contribution of mass from paths with one kink caused by a single intersection with the caustic before adhesion occurs. In this case we see that the shock causes a compression or decompression of mass. This comes from a generalised Ito formula for non $C^2$ functions.

Let $x_0=x_0(t,u,x_0^1)$ where $ x_0^1 \in k^1(u,t)\subset  \mathbb{R}$  and $x_0\in K^1(u,t)$ the curvilinear open set,
\[K^1(u,t) = \left(\Phi_u^{-1} C^{\mathrm{Cool}}_u\right)^{\tilde{}}\cap\left(\left(\lim\limits_n \Phi_t^n\right)^{-1} M_t^{\mathrm{Cool}} \right)^{\tilde{}},\]
\[K^1(t) = \bigcup_{u\in(0,t)} K^1(u,t). \]
If we then integrate over $K^1(t)$, and use the generalised Ito formula for a discontinuous function due to Elworthy, Truman and Zhao \citep{ETZ2, MR2238942},  we get,
\[m^1(0,T) = \int_0^T\hf[ \di t] \int_0^t\di u \int_{k^1(u,t)}\di x_0^1\,T_0^2(x_0) \left|\frac{\partial x_0(t,u,x_0^1)}{\partial (t,u,x_0^1)}\right| \times\]
\[\qquad\qquad\qquad\qquad\qquad \exp\left\{-[\nabla_n\mathcal{S}_u(X(u^+))-\nabla_n\mathcal{S}_u(X(u^{-}))]\right\},\]
where $\mathcal{S}_u(X(u^{\pm}))$ is the action evaluated at the  appropriate minimising pre-image $x_0(i)(x,t)$ just above or below the caustic $C_u$ and 
$\nabla_n \mathcal{S}$ is the normal derivative. The final factor in the integrand is the compression/decompression term coming from the adhesion process. In one dimension this result can be cast as a theorem. In higher dimensions it is only so far a conjecture.

Clearly, the mass involved in the formation of galaxies in the early universe would be given by a sum of terms involving $m^r(0,T)$ for $r=0,1,2,\ldots$. As discussed at the start of this section, the complete adhesion of all matter to the cool part of the Maxwell set precisely destroys the rotation discussed in Section 7. However, if the adhesion were only partial, this would not totally remove the rotation, and could explain the formation of spiral galaxies in the early universe. We will discuss such properties and detailed examples in a forthcoming work on the Burgers equation with voriticity \citep{Vortex}.
\section*{Acknowledgements}
It is a pleasure for AN to thank the Welsh Institute for Mathematical and Computational Sciences (WIMCS) for their financial support in this research.

\def\cprime{$'$}

\end{document}